\documentclass[11pt]{amsart}
\usepackage{amsmath,amssymb,amsthm,graphicx}
\newcommand{\R}{\mathbb{R}}
\newcommand{\Q}{\mathbb{Q}}
\newcommand{\Z}{\mathbb{Z}}

\newtheorem{theorem}{Theorem}
\newtheorem{corollary}[theorem]{Corollary}
\newtheorem*{levy}{Theorem (Levy)}
\theoremstyle{definition}
\newtheorem{claim}{Claim}
\begin{document}
\title[Sharp bounds for harmonic polynomials]{Sharp bounds for the valence of certain harmonic polynomials}
\author{Lukas Geyer}
\address{Lukas Geyer\\Montana State University\\Department of
Mathematics\\ P.O.~Box 172400\\Bozeman, MT 59717--2400\\ USA}
\email{geyer@math.montana.edu}
\subjclass[2000]{26C10; 30C10, 37F10}
\begin{abstract}
  In \cite{KS} it was proved that harmonic polynomials
  $z-\overline{p(z)}$, where $p$ is a holomorphic polynomial of degree
$n > 1$, have at most $3n-2$ complex
  zeros. We show that this bound is sharp for all $n$ by proving a
conjecture of Sarason and Crofoot about the existence of certain extremal
polynomials $p$. We also count the number of equivalence classes of
these polynomials.
\end{abstract}
\maketitle
\section{Introduction and Results} In \cite{KS} it was proved that
harmonic polynomials $z-\overline{p(z)}$, $\deg p = n > 1$, have at
most $3n-2$ complex zeros. An example was provided that this bound is
sharp for $n=3$ and it was conjectured that it is sharp for every
$n$. The main result of this note is the following theorem, which was
conjectured by Crofoot and Sarason \cite[Conj.~1]{KS}.

\begin{theorem}\label{MainThm}
  For every $n>1$ there exists a complex analytic polynomial $p$ of
  degree $n$ and mutually distinct points $z_1, z_2, \ldots, z_{n-1}$ with
  $p'(z_j)=0$ and $\overline{p(z_j)} = z_j$.
\end{theorem}

In particular, this immediately
implies the sharpness of the bound $3n-2$ for every $n > 1$.

\begin{corollary}\label{MainCor}
  For every $n>1$ there exists a complex analytic polynomial $p$ of
degree $n$ such that $\overline{p(z)}-z$ has $3n-2$ zeros.
\end{corollary}

It was proved in \cite{KS} that Theorem \ref{MainThm} implies
Corollary \ref{MainCor}. However, for the convenience of the reader we
will provide an alternative argument to derive this Corollary after
the proof of the Theorem. The main ingredient in the proof of Theorem
\ref{MainThm} is an explicit construction of a ``topological
polynomial'' and a result of Levy that certain topological polynomials
are Thurston-equivalent to polynomials.

In addition we are able to calculate the number of real polynomials
$p$ satisfying the conditions in Theorem \ref{MainThm}, up to equivalence. Two
polynomials $p$ and $q$ are \emph{conjugate} if there is an affine
linear map $T$ with $p = T^{-1} \circ q \circ T$. They are
\emph{equivalent} if there are affine linear maps $S$ and $T$ with
$p=S \circ q \circ T$.

\begin{theorem}\label{CountingThm}
  The number $E_n$ of equivalence classes of real polynomials $p$ of
  degree $n$ having $n-1$ distinct critical points
  $c_1,\ldots,c_{n-1}$ and satisfying $p(c_j) = \overline{c_j}$ for
  $j=1,\ldots,n-1$ is $E_n = C_{\lfloor (n-1)/2 \rfloor}$, where $C_m
  = \frac{1}{m+1} {{2m} \choose {m}}$ is the $m$-th Catalan
number. The number of conjugacy classes is $Q_n = E_n + C_k$ if
$n=4k+3$, and $Q_n = E_n$ if $n \not\equiv 3 \pmod 4$.
\end{theorem}

{\sc Remark.} There are non-real polynomials satisfying the same
equations. In particular, for any real solution $p$, and $\omega$ on
the unit circle, the function $q(z) = \omega p(\omega z)$ is also a
solution. We conjecture that all non-real solutions are equivalent to
real solutions.

The proof of Theorem \ref{CountingThm} relies on Poirier's results
about Hubbard trees for post-critically finite polynomials from
\cite{Po1} and \cite{Po2}.

If one allows rational functions $p$ instead of polynomials, the
number of zeros of $\overline{p(z)}-z$ is at most $5n-5$ (see
\cite{KN}, which also outlines the connection of this problem to
gravitational lensing), and the proof that this bound is sharp is much
easier and very explicit \cite{Rh}.

{\sc Acknowledgment.} I would like to thank Dmitry Khavinson for
helpful discussions.

\section{Proof of Theorem \ref{MainThm}}
We will construct a polynomial with real
coefficients. There is a minor difference between the cases where $n$
is even and the one where $n$ is odd. In the case of even $n$, one of
the points will be real, in the odd case all $z_j$ will be non-real.
Note also that the polynomial $p$ itself will map all $z_j$ to their
complex conjugates, so the second iterate fixes all $z_j$. These
polynomials are special cases of post-critically finite rational
maps. The main ingredient in the proof is Thurston's topological
characterization of rational maps \cite{DH}, and a result of Levy on
topological polynomials. For background see also \cite{BFH} and
\cite{Pi}.

A \emph{topological polynomial} is an orientation-preserving branched
covering $f$ of the sphere to itself with $f^{-1}(\infty) = \{ \infty
\}$. Critical points are points where $f$ is not locally injective,
the \emph{post-critical set} $P(f)$ is the closure of $\{f^n(c):
n>1,\, c \text{ critical}\}$. An orientation-preserving branched
covering of the sphere is \emph{post-critically finite} if the
post-critical set contains finitely many points, or equivalently if
all critical points are eventually periodic under iteration. Two
post-critically finite maps $f$ and $g$ are \emph{Thurston equivalent}
if there exist orientation-preserving homeomorphisms $\phi$ and
$\psi$, homotopic rel $P(f)$, with $\phi \circ f = g \circ \psi$.

Thurston gave a criterion for post-critically finite branched
coverings to be equivalent to rational maps, and Levy proved the
following about topological polynomials (see \cite{Le},
\cite[Cor.~5.13]{BFH}).
\begin{levy}
 If $f$ is a topological polynomial such that every critical point
eventually lands in a periodic cycle containing a critical point, then
$f$ is Thurston-equivalent to a polynomial. The polynomial is unique
up to affine conjugation.
\end{levy}

We will first construct a topological polynomial $f$ satisfying
$f(\overline{z}) = \overline{f(z)}$ for all $z$, and $f(c_j) =
\overline{c_j}$ for all critical points $c_j$. The cases of even and
odd degrees are slightly different, so we will start with the case of
odd degree $n=2m+1$ for some $m \geq 1$. Let $g$ be a polynomial of
degree $2m+1$ with real coefficients, $m$ distinct critical points
$c_1,\ldots, c_m$ in the upper halfplane such that the corresponding
critical values $v_j = g(c_j)$ are mutually distinct points in the
lower halfplane. E.g., one could choose a small constant $\epsilon >
0$, define $c_{2k-1} = ik$, $c_{2k} = i(k+\epsilon)$, and
\[
g(z) = - \int_0^z \prod_{k=1}^m (\zeta^2 - c_k^2) d\zeta.
\]
It is easy to check that $g$ satisfies all the conditions if
$\epsilon$ is sufficiently small. Now choose a homeomorphism $h$ of
the closed lower halfplane with $h(v_j) = \overline{c_j}$. Extend $h$
by reflection to a homeomorphism of the plane and define $f = h \circ
g$. Then $f$ is a topological polynomial of degree $2m+1$, commuting
with reflection in the real line, mapping all its critical points to
their complex conjugates. The post-critical set of $f$ consists of all
the critical points
$c_1,\ldots,c_m,\overline{c_1},\ldots,\overline{c_m}$ and $\infty$. By
Levy's result $f$ is Thurston-equivalent to a polynomial $p$. The
uniqueness part of Levy's theorem implies that $p$ is conjugate to
$\overline{p}$, so $p$ is conjugate to a real polynomial, and we may
as well assume that $p$ itself is real. It also follows that both
$\phi$ and $\psi$ are symmetric with respect to the real line. By the
definition of Thurston-equivalence we know that $z_j = \phi(c_j) =
\psi(c_j)$, so $p$ has critical points $z_1,\ldots,z_m,
\overline{z_1}, \ldots, \overline{z_m}$, and $p(z_j) =
\overline{z_j}$ for all $j$.

If the degree $n$ is even, we have to modify the initial step of the
construction. Write $n=2m+2$ with some $m \geq 0$, and find a real
polynomial $g$ of degree $n$ with distinct simple critical points $c_0
\in \R$, $c_1,\ldots,c_m$ in the upper halfplane, such that $v_j = g(c_j)$
are all distinct and in the lower halfplane for $1 \leq j \leq m$. Here one
can take
\[
g(z) = \int_0^z (\zeta-1) \prod_{k=1}^m (\zeta^2 - c_k^2) d\zeta
\]
with $c_k$ chosen as above for $k \geq 1$. Now choose a homeomorphism $h$
of the closed lower halfplane onto itself such that $h(g(c_0)) = c_0$
and $h(g(c_k)) = \overline{c_k}$ for $k \geq 1$. The rest of the proof is
identical to the case of odd degree.\qed

{\sc Proof of Corollary \ref{MainCor}.} We will show that
$z-\overline{p(z)}$ has at least $3n-2$ zeros. Together with the result
from \cite{KS}, this shows that it has exactly $3n-2$ zeros.  First we
observe that $p$ has no degenerate fixed points, i.e.~all fixed points
of $z\mapsto \overline{p(z)}$ satisfy $|p'(z)| \neq 1$. (For
holomorphic maps, a fixed point is degenerate if the derivative is
exactly 1, for anti-holomorphic maps it is degenerate whenever the
derivative has modulus 1.) If this were the case, then
$\overline{p(\overline{p(z)})}$ would have a parabolic fixed point,
which would have a critical point in its attracting basin. However,
all critical points of this maps are mapped to the critical points of
$p$ and those are fixed. Having only non-degenerate fixed points, the
sum of the Lefschetz indices of the fixed points in the Riemann sphere
has to be equal to $-n+1$, as the induced maps on 0- and 2-dimensional
homology are the identity and multiplication by $-n$, respectively. It
already has super-attracting fixed points at $\infty$ and $z_1,
\ldots, z_{n-1}$. The Lefschetz index of a super-attracting fixed point is $1$,
thus the sum of the indices of all those fixed points is $n$. We
conclude that $p$ must have at least $2n-1$ fixed points with index
$-1$, i.e. repelling fixed points for $\overline{p}$. This proves
sharpness. \qed

\section{Counting Conjugacy and Equivalence Classes}
We will use results about Hubbard trees by Poirier (see \cite{Po1} and
\cite{Po2}) in order prove Theorem \ref{CountingThm}. Hubbard
trees were first introduced in \cite{BFH} in order to classify
post-critically finite polynomials without periodic critical points.
In that case the Julia set is a dendrite, and the Hubbard tree is the
embedded subtree of the Julia set spanned by the post-critical set. It
has a finite number of vertices and edges.  The polynomial induces a
map on the tree, and conversely an embedded finite tree with a map
satisfying certain conditions gives rise to a post-critically finite
polynomial, unique up to affine conjugation.  Poirier generalized this
to include the case of post-critically finite polynomials with
critical periodic points. In this case all Fatou components are simply
connected with a distinguished base point, the unique preperiodic
point in the component. The Hubbard tree in this case intersects both
the Julia and the Fatou set, and the part in a Fatou set consists of
preimages of radial lines under the Riemann map mapping the component
to the unit disk and the basepoint to 0. This construction determines
an angle between any two edges meeting at a point in the Fatou set.
The important information needed to reconstruct the map from the
(marked) tree is the local degree at every vertex, the angles between
any two edges meeting at a Fatou vertex, and the mapping on
vertices. There is no well-defined angle between edges meeting at
points in the Julia set, but there is a well-defined cyclic order, and
we define the ``angles'' between the edges to be all the same. In that
way we end up with an angle function on the tree defined for any pair
of incident edges.

An abstract Hubbard tree is a finite tree embedded in the plane,
together with a mapping $F:V\to V$ on the vertices, a local degree
function $\delta:V\to\{ 1,2,\ldots \}$ on the vertices, and an angle function $\alpha$
defined for any pair of incident edges $l_k, l_j$. It has to satisfy
$\alpha(l_k,l_j)\in \R\setminus \Q$, $\alpha(l_j,l_k) = -\alpha(l_k,l_j)$, and $\alpha(l_j, l_k) +
\alpha(l_k, l_m) = \alpha(l_j,l_m)$ whenever $l_k, l_j, l_m$ meet at a vertex.
The map $F$ satisfies $F(v) \neq F(w)$ whenever $v$ and $w$ are adjacent
vertices. In this way $F$ induces a mapping on the edges, mapping the
edge $[v,w]$ homeomorphically to the unique path between $F(v)$ and
$F(w)$, i.e.~mapping an edge to a union of edges. If two edges $l_j$
and $l_k$ meet at $v$, then $F(l_j)$ and $F(l_k)$ meet at $F(v)$, and
we require that $\alpha(F(l_j),F(l_k)) = \delta(v) \alpha(l_j,l_k)$, which is just
the fact that holomorphic maps of local degree $\delta$ multiply angles by
$\delta$.  A vertex is a Fatou vertex if it eventually lands in a periodic
vertex $v$ with $\delta(v)>1$.  All other vertices are Julia vertices. The
last and not quite as obvious requirement for the Hubbard tree is that
$F$ is expanding, i.e.~whenever $v$ and $w$ are adjacent Julia
vertices, there exists $n$ such that $F^n(v)$ and $F^n(w)$ are not
adjacent. 

The \emph{degree} of the Hubbard tree is $d = 1 + \sum_{v \in V} (\delta(v)-1)$.
Poirier's main result in \cite{Po2} is that any such abstract Hubbard
tree corresponds to a polynomial of degree $d$, unique up to affine
conjugation. The uniqueness part easily implies that any symmetries of
the Hubbard tree correspond to symmetries of the polynomial.

{\sc Proof of Theorem \ref{CountingThm}.} We will first prove the
result about the number of conjugacy classes. By Poirier's theorem we
only have to count the number of corresponding Hubbard trees up to
symmetry. We will first establish some general properties of the
corresponding Hubbard trees, after which we will treat odd and even
degrees separately. The following claim gives a complete
characterization of the Hubbard trees in question, the rest of the
argument is counting.

\begin{claim}
The induced map on the Hubbard tree is $z \mapsto \overline{z}$. The
only vertices in the Hubbard tree are the critical points. All angles
are $1/3$ or $2/3$. Conversely, any finite tree $T$ which is symmetric with
respect to complex conjugation, has local degree $2$ at every vertex,
with the map $f:T\to T$, $f(z) = \overline{z}$, arises as a Hubbard tree
of a polynomial satisfying the requirements in Theorem \ref{CountingThm}.
\end{claim}

{\sc Proof.} Let $c_1$ and $c_2$ be critical points such that the path
$\gamma$ between them does not contain any other critical point. Then
$\gamma$ will be mapped homeomorphically to the path connecting
$p(c_1) = \overline{c_1}$ and $p(c_2) = \overline{c_2}$, which is the
complex conjugate of $\gamma$. Since the whole post-critical set
consists of critical points, this argument covers the whole Hubbard
tree. Any non-critical vertices would have to be branch points, and
the induced map $z \mapsto \overline{z}$ reverses the cyclic
order of edges at branch points, contradicting the fact that $p$ is
locally conformal at all non-critical vertices. Finally, let $l_1$ and
$l_2$ be two edges meeting at a critical vertex $c$. Then $\delta(c) =
2$, and $2 \alpha(l_1,l_2) = \alpha(p(l_1),p(l_2)) =
\alpha(\overline{l_1}, \overline{l_2}) = -\alpha(l_1,l_2)$ in $\R /
\Z$, thus $3 \alpha(l_1, l_2) \in \Z$. 

Conversely, if all local degrees are 2, and all angles are $1/3$ or
$2/3$, complex conjugation $f(z) = \overline{z}$ satisfies $\delta(v)
\alpha(l_1,l_2) = \alpha(f(l_1),f(l_2))$ for any two edges meeting at
a vertex $v$. The property that $f$ on the tree $T$ be expanding is
trivially satisfied as there are no Julia vertices at
all. $\R$-symmetric Hubbard trees give rise to $\R$-symmetric
polynomials, so this shows the last part of the claim.

\begin{figure}
\includegraphics*[height=2in]{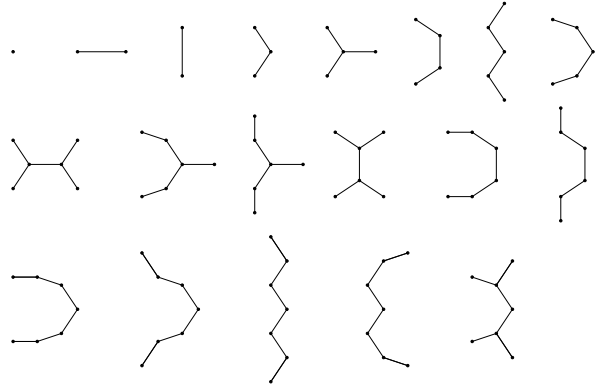}
\caption{Hubbard trees for degrees $2 \le n \le 8$. The real line is
the horizontal line of symmetry for each tree, and the induced mapping
on the trees is complex conjugation. All angles are multiples of $2\pi/3$.}
\label{FigTrees}
\end{figure}

\begin{claim}
The number of vertices on the real line is at most 2. There is a
bijective correspondence between Hubbard trees with 0 real vertices
and those with 2 real vertices.
\end{claim}
{\sc Remark.} The bijective correspondence replaces the vertical
edge which intersects the real line by a horizontal edge which is a
real interval. See Figure \ref{FigTrees} for examples of Hubbard trees
up to degree 8.

{\sc Proof.} Let $c_1$ and $c_2$ be real critical points. Since the
Hubbard tree is a tree symmetric with respect to the real line, the
subtree spanned by $c_1$ and $c_2$ has to be the real interval between
them. Since no angle in the Hubbard tree is $1/2$, there can be no
other vertex in that interval, i.e.~there can not be other critical
points between $c_1$ and $c_2$. This shows the first claim. In order
to construct the bijective correspondence, we only need to give the
modification of the tree in the upper halfplane. The part in the lower
halfplane is determined by symmetry and the map on the tree is always
complex conjugation. If there are two real critical points $c_1<c_2$,
let $T_1$ and $T_2$ be the subtrees in the upper halfplane attached at
$c_1$ and $c_2$, resp. The correspondence is given by replacing the
edge $l=[c_1, c_2]$ by $\tilde{l}=[-i,i]$, and attaching $T_1$ and
$T_2$ to the left and right of $i$, resp., so that $\alpha(T_2, T_1) =
\alpha(T_1,\tilde{l}) = \alpha(\tilde(l),T_2) = 1/3$.  The inverse
of this operation is given by replacing the edge which crosses the
real line by a real interval and attaching the left and right subtrees
at the endpoints of this interval, again with angles $1/3$. Observe
also that any symmetry with respect to the imaginary axis will be
preserved by this correspondence and its inverse.

It remains to count the Hubbard trees. Let us start with the easier
case of even degree $n=2m+2$ with $m \ge 0$. Here we have an odd
number of critical points, so by Claim 2 there is exactly one real
critical point $c_1$. The number of edges meeting at this point must
be even by symmetry. In the case $m=0$, there is only one critical
point, so there are no edges. In all other cases the number of edges
meeting at $c_0$ must be $2$, since it can be at most 3 by Claim
1. Also the angles between the real axis and the edges must be $1/6$
and $1/3$, respectively. There are exactly two choices, and they
correspond to each other under the symmetry $z \mapsto -z$. Since we
count Hubbard trees up to symmetry, we may fix an orientation. It also
follows that all those Hubbard trees have no rotational symmetry. We
can view the part of the tree in the upper halfplane as a finite part
of an infinite rooted binary tree. The number of vertices, excluding
the root, is $m$, and it is well-known (and easy to proof from the
recursive formula $C_m = \sum\limits_{k=0}^{m-1} C_k C_{m-1-k}$) that there
are $C_m = \frac1{m+1} { {2m} \choose m}$ of those.

In the case of odd degree $n=2m+1$ there are $2m$ critical points. Let
us count the number of trees which have no vertex on the real
line. Then the part in the upper halfplane is again a finite part of a
rooted binary tree with $m$ vertices, so there are $C_m$ of
those. However, for every such tree the image under $z \mapsto -z$ has
the same form. If $m=2k$, this implies that there are $C_m/2$ up to
symmetry. However, if $m=2k+1$, there are trees which are symmetric
under $z \mapsto -z$. The number of those corresponds to the number of
subtrees with $k$ vertices attached to the right of the first vertex
in the upper halfplane. This means that there are $C_k$ symmetric
trees, so we have $(C_m+C_k)/2$ inequivalent Hubbard trees in this
case. By Claim 2, we have the same number of trees with 2 vertices on
the real line, so the number of inequivalent Hubbard trees is
$C_{2k}$ for $n=4k+1$ and $C_{2k+1} + C_k$ for $n=4k+3$.

Let us now count the number of equivalence classes. The equivalence of
two solutions $p$ and $q$ implies the existence of affine linear
maps $S$ and $T$ with $ S p S^{-1} = S q $, i.e.~the existence of an
affine map $R=ST$ such that $R$ maps the critical points of $q$ to
itself. In particular this implies that $RqR^{-1}$ is again a
solution. The only rotations $R$ which map one of our symmetric
Hubbard trees to a symmetric Hubbard tree are $R(z)=i^p z$ with $p \in
\{1,2,3\}$. In the case $R(z)=-z$, the map $p=Rq$ has the same
underlying Hubbard tree as $q$ with the map $z \mapsto -\overline{z}$
on it. This corresponds to the Hubbard tree of $q$ rotated by 90
degrees with complex conjugation as the map. This can only happen if
$n=4k+3$, and the number of trees which admit this is
$C_k$. The case $R(z) = \pm i$ can actually not occur, since the map
$p=Rq$ would have the same Hubbard tree and the map on it would be
reflection in the line $(1 \pm i) \R$. However, none of the Hubbard
trees has one of those lines as a line of symmetry.

Combining all this, the number of equivalence classes is the same as
the number of conjugacy classes except in the case $n=4k+3$, where it
becomes $C_{2k+1}$. Thus the number of equivalence classes is
$Q_n = C_{\lfloor (n-1)/2 \rfloor}$ for all $n$.
\qed

\end{document}